\documentclass[a4paper]{article}
\usepackage[english]{babel}
\usepackage[utf8]{inputenc}
\usepackage{amsmath}
\usepackage{graphicx}
\usepackage[colorinlistoftodos]{todonotes}

\usepackage{amsthm}
\usepackage{amsfonts}
\usepackage{amssymb}
\usepackage{amscd}
\usepackage{stmaryrd}
\usepackage{float}
\usepackage{relsize}
\usepackage{geometry}
\usepackage[font=small,justification=centering]{caption}

\newcommand{\A}{\ensuremath \mathcal{A}}
\newcommand{\B}{\ensuremath \mathcal{B}}
\newcommand{\Q}{\ensuremath \mathbb{Q}}

\newcommand{\C}{\ensuremath \mathbb{C}}
\newcommand{\Z}{\ensuremath \mathbb{Z}}
\newcommand{\N}{\ensuremath \mathbb{N}}
\newcommand{\isom}{\ensuremath \cong}

\title{Transposing Noninvertible Polynomials}

\author{Nathan Cordner}

\date{\today}

\begin{document}
\newgeometry{top = 1 in, left = 1 in, right = 1 in, bottom = 1 in}

\maketitle

\begin{abstract}
Landau-Ginzburg mirror symmetry predicts isomorphisms between graded Frobenius algebras (denoted $\A$ and $\B$) that are constructed from a nondegenerate quasihomogeneous polynomial $W$ and a related group of symmetries $G$. Duality between $\A$ and $\B$ models has been conjectured for particular choices of $W$ and $G$. These conjectures have been proven in many instances where $W$ is restricted to having the same number of monomials as variables (called \emph{invertible}). Some conjectures have been made regarding isomorphisms between $\A$ and $\B$ models when $W$ is allowed to have more monomials than variables. In this paper we show these conjectures are false; that is, the conjectured isomorphisms do not exist. Insight into this problem will not only generate new results for Landau-Ginzburg mirror symmetry, but will also be interesting from a purely algebraic standpoint as a result about groups acting on graded algebras.
\end{abstract}

\section{Introduction}

\noindent Physicists conjectured some time ago that that to each quasihomogeneous (weighted homogeneous) polynomial $W$ with an isolated singularity at the origin, and to each admissible group of symmetries $G$ of $W$, there should exist two different physical ``theories," (called the Landau-Ginzburg $\A$ and $\B$ models, respectively) consisting of graded Frobenius algebras (algebras with a nondegenerate pairing that is compatible with the multiplication).   The $\B$-model theories have been constructed  \cite{IV, Ka1, Ka2, Ka3, Kra09} and correspond to an ``orbifolded Milnor ring."  The $\A$-model theories have also been constructed \cite{FJR07} and are a special case of what is often called ``FJRW theory."  We will not address these in this paper, but in many cases, these theories can be extended to whole families of Frobenius algebras, called \emph{Frobenius manifolds}.
\vspace{0.1 in}

\noindent For a large class of these polynomials (called \emph{invertible}) Berglund-H\"ubsch \cite{BH}, Henningson \cite{Hen}, and Krawitz  \cite{Kra09} described the construction of a dual (or transpose) polynomial $W^T$ and a dual group $G^T$.  The Landau-Ginzburg mirror symmetry conjecture states that the $\A$-model of a pair $W,G$ should be isomorphic to the $\B$-model of the dual pair $W^T,G^T$.  This conjecture has been proved in many cases in papers such as \cite{Kra09} and \cite{FJJS11}, although the proof of the full conjecture remains open.
\vspace{0.1 in}

\noindent It has been further conjectured that the Berglund-H\"ubsch-Henningson-Krawitz duality transform should extend to large classes of noninvertible polynomials and that Landau-Ginzburg mirror symmetry should also hold for these polynomials.  In this paper we investigate some candidate mirror pairs of noninvertible polynomials and show that many obvious candidates for mirror duality cannot satisfy mirror symmetry.

\vspace{0.1 in}

\noindent To approach this problem, we study the $\A$ and $\B$ models as graded vector spaces and inspect how the symmetry groups act on these spaces. Insight into this problem will not only generate new results for Landau-Ginzburg mirror symmetry, but will also be interesting from a purely algebraic standpoint as a result about groups acting on graded algebras.
\vspace{0.1 in}

\noindent One case of mirror symmetry that has been verified for all invertible polynomials is when the $\A$-model is constructed from an invertible polynomial $W$ with its maximal group of symmetries and the $\B$-model is constructed from the corresponding transpose polynomial with the trivial group of symmetries. This is sometimes denoted $\A_{W, G_{W}^{max}} \isom \B_{W^{T}, \{0\}}$. This intuition stemming from invertible polynomials motivated two conjectures about isomorphisms between $\A$ and $\B$ models built from noninvertible polynomials. We often refer to polynomials for which the $\A$ and $\B$ models exist as \emph{admissible}.
\vspace{0.1 in}

\noindent \textbf{Conjecture 1}. For any admissible (not necessarily invertible) polynomial $W$ in $n$ variables, there exists a corresponding admissible polynomial $W^{T}$ in $n$ variables satisfying $\A_{W,G_{W}^{max}} \isom \B_{W^{T}, \{0\}}$.
\vspace{0.1 in}

\noindent Note that this conjecture includes the collection of noninvertible polynomials, which are allowed to have more monomials than variables. In Section 3.1 we show that this conjecture is false. By relaxing the restriction on the number of variables that $W^{T}$ is allowed to have, we obtain a second conjecture.

\vspace{0.1 in}

\noindent \textbf{Conjecture 2}. For any admissible $W$, there is a corresponding admissible $W^{T}$ satisfying $\A_{W,G_{W}^{max}} \isom \B_{W^{T}, \{0\}}$.

\vspace{0.1 in}

\noindent In Section 3.2 we look at an example of a particular noninvertible polynomial, and expand our search space for finding a suitable $W^{T}$. We develop some formulas and show that they rule out the existence of $W^{T}$ in a few more cases that were not considered in Conjecture 1. Thereby we also establish that Conjecture 2 is unlikely to be true in general.

\section{Preliminaries}

Here we will introduce some of the concepts needed to explain the theory of this paper. 

\subsection{Admissible Polynomials}

\noindent \textbf{Definition}. For a polynomial $W \in \C[x_{1},\dots,x_{n}]$, we say that $W$ is \textit{nondegenerate} if it has an isolated critical point at the origin.

\vspace{0.1 in}

\noindent \textbf{Definition}. Let $W \in \C[x_{1},\dots,x_{n}]$. We say that $W$ is \textit{quasihomogeneous} if there exist positive rational numbers $q_{1},\dots,q_{n}$ such that for any $c \in \C$, $W(c^{q_{1}}x_{1},\dots,c^{q_{n}}x_{n}) = c W(x_{1},\dots,x_{n})$.

\vspace{0.1 in}

\noindent We often refer to the $q_{i}$ as the \textit{quasihomogeneous weights} of a polynomial $W$, or just simply the \textit{weights} of $W$, and we write the weights in vector form $J = (q_{1}, \dots, q_{n})$. 

\vspace{0.1 in}

\noindent \textbf{Definition}. $W \in \C[x_{1},\dots,x_{n}]$ is \textit{admissible} if $W$ is both nondegenerate and quasihomogeneous, with the weights of $W$ being unique.

\vspace{0.1 in}

\noindent We will use the following result about admissible polynomials later in the paper. 
\vspace{0.1 in}

\noindent \textbf{Proposition 2.1.6 of \cite{FJR07}}. If $W \in \C[x_{1},\dots,x_{n}]$ is admissible, and contains no monomials of the form $x_{i}x_{j}$ for $i \ne j$, then the $q_{i}$ are bounded above by $\frac{1}{2}$. 
\vspace{0.1 in}

\noindent Because the construction of $\A_{W,G}$ requires an admissible polynomial, we will only be concerned with admissible polynomials in this paper. In order for a polynomial to be admissible, it needs to have at least as many monomials as variables. Otherwise its quasihomogeneous weights cannot be uniquely determined. We now state the main subdivision of the admissible polynomials.
\vspace{0.1 in}

\noindent \textbf{Definition}. Let $W$ be an admissible polynomial. We say that $W$ is \textit{invertible} if it has the same number of monomials as variables. If $W$ has more monomials than variables, then it is \textit{noninvertible}.

\vspace{0.1 in}

\noindent Admissible polynomials with the same number of variables as monomials are called invertible since their associated exponent matrices (which we define in the next section) are square and invertible.

\subsection{Dual Polynomials}

\noindent We will now introduce the idea of the transpose operation for invertible polynomials.
\vspace{0.1 in}

\noindent \textbf{Definition}. Let $W \in \C[x_{1},\dots, x_{n}]$. If we write $W = \sum_{i = 1}^{m} c_{i} \prod_{j = 1}^{n} x_{j}^{a_{ij}}$, then the associated \textit{exponent matrix} is defined to be $A = (a_{ij})$. 
\vspace{0.1 in}

\noindent From this definition we notice that $n$ is the number of variables in $W$, and $m$ is the number of monomials in $W$. $A$ is an $m \times n$ matrix. Thus when $W$ is invertible, we have that $m = n$ which implies that $A$ is square. One can show, without much work, that this square matrix is invertible if the polynomial $W$ is quasihomogeneous with unique weights. When $W$ is noninvertible, $m > n$. $A$ then has more rows than columns. 

\vspace{0.1 in}

\noindent Observe that if a polynomial is invertible, then we may rescale all nonzero coefficients to 1. So there is effectively a one-to-one correspondence between exponent matrices of invertible polynomials and the polynomials themselves.

\vspace{0.1 in}

\noindent \textbf{Definition}. Let $W$ be an invertible polynomial. If $A$ is the exponent matrix of $W$, then we define the \textit{transpose polynomial} to be the polynomial $W^{T}$ resulting from $A^{T}$. By the classification in \cite{KS}, $W^{T}$ is again a nondegenerate, invertible polynomial.
\vspace{0.1 in}

\noindent We now have reached our fundamental problem. When a polynomial $W$ is noninvertible, its exponent matrix $A$ is no longer square. Taking $A^{T}$ yields a polynomial with fewer monomials than variables, which is not admissible. Therefore, we will require a different approach to define what the transpose polynomial should be for noninvertibles.

\subsection{Symmetry Groups and Their Duals}

\noindent \textbf{Definition}. Let $W$ be an admissible polynomial. We define the \textit{maximal Abelian symmetry group} of $W$ to be $G_{W}^{max} = \{(\zeta_{1}, \dots, \zeta_{n}) \in (\C^\times)^{n} \mid W(\zeta_{1}x_{1},\dots,\zeta_{n}x_{n}) = W(x_{1},\dots,x_{n})    \}$.

\vspace{0.1 in}

\noindent The proofs of Lemma 2.1.8 in \cite{FJR07} and Lemma 1 in \cite{ABS11} observe that $G_{W}^{max}$ is finite and that each coordinate of every group element is a root of unity. The group operation $\circ$ in $G_{W}^{max}$ is coordinate-wise multiplication. That is,
\begin{align*}
(e^{2\pi i \theta_{1}},\dots,e^{2\pi i \theta_{n}}) \circ (e^{2\pi i \phi_{1}},\dots,e^{2\pi i \phi_{n}}) = (e^{2\pi i (\theta_{1} + \phi_{1})},\dots,e^{2\pi i (\theta_{n} + \phi_{n})}).
\end{align*}
Equivalently, in additive notation we can write $(\theta_{1}, \dots, \theta_{n}) + (\phi_{1}, \dots, \phi_{n}) = (\theta_{1} + \phi_{1}, \dots, \theta_{n} + \phi_{n}) \mod \Z$. The map $(e^{2\pi i \theta_{1}},\dots,e^{2\pi i \theta_{n}}) \mapsto (\theta_{1}, \dots, \theta_{n}) \mod \Z$ gives a group isomorphism. Using additive notation, we will often write $G_{W}^{max} = \{ g \in (\Q/\Z)^{n} \mid Ag \in \Z^{m} \}$, where $A$ is the $m \times n$ exponent matrix of $W$.
\vspace{0.1 in}

\noindent \textbf{Definition}. In this notation, $G_{W}^{max}$ is a subgroup of $(\Q/\Z)^{n}$ with respect to coordinate-wise addition. For $g \in G_{W}^{max}$, we write $g = (g_{1}, \dots, g_{n})$ where each $g_{i}$ is a rational number in the interval [0,1). The $g_{i}$ are called the \textit{phases} of $g$.

\vspace{0.1 in}

\noindent The following definition of the transpose group is due to Krawitz and Henningson \cite{Kra09, Hen}. 

\vspace{0.1 in}

\noindent \textbf{Definition}. Let $W$ be an invertible polynomial, and let $A$ be its associated exponent matrix. The \textit{transpose group} of a subgroup $G \le G_{W}^{max}$ is the set $G^{T} = \{g \in G_{W^{T}}^{max} \mid g A h^{T} \in \Z \text{ for all } h \in G   \}$.
\vspace{0.1 in}

\noindent Since this relies on knowing what $W^{T}$ is, this definition currently does not extend to noninvertible polynomials. The following is a list of common results for the transpose group. 
\vspace{0.1 in}

\noindent \textbf{Proposition 2 of \cite{ABS11}}. Let $W$ be an invertible polynomial with weights vector $J$, and let $G \le G_{W}^{max}$. 

(1) $(G^{T})^{T} = G$,

(2) $\{0\}^{T} = G_{W^{T}}^{max}$ and $(G_{W}^{max})^{T} = \{0\}$,

(3) $\langle J \rangle^{T} = G_{W^{T}}^{max} \cap \text{SL}(n,\C)$ where $n$ is the number of variables in $W$,

(4) if $G_{1} \le G_{2}$, then $G_{2}^{T} \le G_{1}^{T}$ and $G_{2} / G_{1} \isom G_{1}^{T} / G_{2}^{T}$.

\subsection{Some Notes on $\A$ and $\B$ Models}

Landau-Ginzburg $\A$ and $\B$ models are algebraic objects that are endowed with many levels of structure. In this paper, we will chiefly be concerned with their structure as graded vector spaces, although we will also occasionally consider their Frobenius algebra structure. For the benefit of the reader, we will give a formal definition of a Frobenius algebra. 

\vspace{0.1 in}

\noindent \textbf{Definition}. An \textit{algebra} is a vector space $A$ over a field of scalars $F$ (in our case it is $\C$), together with a multiplication $\cdot : A \times A \rightarrow A$ that satisfies for all $x,y,z \in A$ and $\alpha,\beta \in F$

$\bullet$ Right distributivity:  $(x + y) \cdot z = x \cdot z + y \cdot z$,

$\bullet$ Left distributivity:  $x \cdot (y + z) = x \cdot y + x \cdot z$,

$\bullet$ Compatability with scalars: $(\alpha x) \cdot (\beta y) = (\alpha\beta) (x \cdot y)$.

\noindent We further require the multiplication to be associative and commutative, and for $A$ to have a unity $e$ such that $e \cdot x = x$ for all $x \in A$. 

\vspace{0.1 in}

\noindent We also define a pairing operation $\langle \cdot, \cdot \rangle: A \times A \rightarrow F$ that is 

$\bullet$ Symmetric:  $\langle x, y \rangle = \langle y, x \rangle$,

$\bullet$ Linear:  $\langle \alpha x + \beta y, z \rangle = \alpha \langle x, z \rangle + \beta \langle y, z \rangle$,

$\bullet$ Nondegenerate:  for every $x \in A$ there exists $y \in A$ such that $\langle x, y \rangle \ne 0$. 

\noindent If the pairing further satisfies the \textit{Frobenius property}, meaning that $\langle x \cdot y, z \rangle = \langle x, y\cdot z \rangle$ for all $x,y,z \in A$, then we call $A$ a \textit{Frobenius algebra}.

\vspace{0.1 in}

\noindent We will only develop the theory needed for the proofs in Section 3. We refer the interested reader to \cite{FJR07} for more details on the construction of the $\A$-model. \cite{FJJS11}, \cite{Kra09}, and \cite{Tay13} also contain more information on constructing $\A$ and $\B$ models, and related isomorphisms. We will start by discussing the $\B$-model. 

\vspace{0.1 in}

\noindent \textbf{Definition}. $\mathcal{Q}_{W} = \C[x_{1},\dots,x_{n}] / (\frac{\partial W}{\partial x_{1}}, \dots, \frac{\partial W}{\partial x_{n}})$ is called the \textit{Milnor ring} of $W$ (or \textit{local algebra} of $W$).

\vspace{0.1 in}

\noindent \textbf{Definition}. We define the \textit{unorbifolded} $\B$-model to be $\B_{W, \{0\}} = \mathcal{Q}_{W}$.

\vspace{0.1 in}

\noindent We will think of the unorbifolded $\B$-model as a graded vector space over $\C$. The degree of a monomial in $\mathcal{Q}_{W}$ is given by $\deg(x_{1}^{a_{1}}x_{2}^{a_{2}}\dots x_{n}^{a_{n}}) = 2\sum_{i}^{n} a_{i}q_{i}$. This defines a grading on the basis of $\mathcal{Q}_{W}$. We note the following:

\vspace{0.1 in}

\noindent \textbf{Theorem 2.6 of \cite{Tay13}}. If $W$ is admissible, then $\mathcal{Q}_{W}$ is finite dimensional.

\vspace{0.1 in}

\noindent We will need two results about the unorbifolded $\B$-model. First, $\dim(\B_{W,\{0\} }) = \prod_{i = 1}^{n} \left( \frac{1}{q_{i}} - 1 \right)$. Second, the highest degree of its graded pieces is $2\sum_{i = 1}^{n} \left(1 - 2q_{i} \right)$. (See Section 2.1 of \cite{Kra09})

\vspace{0.1 in}

\noindent We will now develop some needed ideas about $\A$-models. 
\vspace{0.1 in}

\noindent \textbf{Definition}. Let $W$ be an admissible polynomial with weights vector $J = (q_{1}, \dots, q_{n})$, and let $G \le G_{W}^{max}$. Then $G$ is \textit{admissible} if $J \in G$. 
\vspace{0.1 in}

\noindent We note that since $W$ is quasihomogeneous, we have that $AJ^{T} = (1,\dots,1)^{T} \in \Z^{m}$. Thus $J \in G_{W}^{max}$.

\vspace{0.1 in}

\noindent The construction of the $\A$-model requires that $G$ be an admissible group. From parts (3) and (4) of the proposition in Section 2.3, the corresponding condition for the $\B$-model is that $G^{T} \le G_{W^{T}}^{max} \cap \text{SL}(n,\C)$.

\vspace{0.1 in}

\noindent \textbf{Definition}. Let $W \in \C[x_{1},\dots,x_{n}]$ be admissible, and let $g = (g_{1}, \dots, g_{n}) \in G_{W}^{max}$. The \textit{fixed locus} of the group element $g$ is the set fix$(g) = \{x_{i} \mid g_{i} = 0  \}$.

\vspace{0.1 in}

\noindent We now state how $G$ acts on the Milnor ring.

\vspace{0.1 in}

\noindent \textbf{Definition}. Let $W$ be an admissible polynomial, and let $g \in G_{W}^{max}$. We define the map $g^{*}: \mathcal{Q}_{W} \rightarrow \mathcal{Q}_{W}$ by $g^{*}(m) = \det(g) m \circ g$. (Here we think of $g$ as being a diagonal map with multiplicative coordinates)

\vspace{0.1 in}

\noindent \textbf{Definition}. Let $W$ be an admissible polynomial, and let $G \le G_{W}^{max}$. Then the \emph{$G$-invariant subspace} of $\mathcal{Q}_{W}$ is defined to be $\mathcal{Q}_{W}^{G} = \{m \in \mathcal{Q}_{W} \mid g^{*}(m) = m \text{ for each } g \in G \}$.

\vspace{0.1 in}

\noindent \textbf{Definition}. Let $W$ be an admissible polynomial, and $G$ an admissible group. We define $\A_{W,G} = \mathlarger{\bigoplus\limits_{g \in G}} \left( \mathcal{Q}_{W|_{\text{fix}(g)}} \right)^{G}$, where $( \cdot )^{G}$ denotes all the $G$-invariants. This is called the $\A$-model \textit{state space}.

\vspace{0.1 in}

\noindent We further note that the state space of the \textit{orbifolded} $\B$-model $\B_{W,G}$ is constructed similarly, but with the condition that $G \le G_{W}^{max} \cap \text{SL}(n,\C)$. If we let $G = \{0\}$, then the formula yields the Milnor ring of $W$ as expected. The grading on the $\A$-model, which will will define in a moment, differs from the $\B$-model grading; but as graded vector spaces, the $\A$ and $\B$ models are very much related. 
\vspace{0.1 in}

\noindent We will not discuss many details of constructing the state space here. For further treatment of this topic, we refer the reader to Section 2.4 of \cite{Tay13}. A brief comment on notation:  we represent basis elements of $\A_{W,G}$ in the form $\lfloor m; g \rceil$, where $m$ is a monomial and $g$ is a group element. 

\vspace{0.1 in}

\noindent \textbf{Definition}. The $\A$-model degree of a basis element $\lfloor m ; g \rceil$ is defined to be $\deg( \lfloor m ; g \rceil ) =  \dim(\text{fix}(g)) + 2 \sum_{i=1}^{n}(g_{i} - q_{i})$, where $g = (g_{1},\dots,g_{n})$ with the $g_{i}$ chosen such that $0 \le g_{i} < 1$ and $J = (q_{1},\dots,q_{n})$ is the vector of quasihomogeneous weights of $W$. (See Section 2.1 of \cite{Kra09})
\vspace{0.1 in}

\noindent Finally, we state one important theorem for $\A$-model isomorphisms.

\vspace{0.1 in}

\noindent \textbf{Theorem in Section 7.1 of \cite{Tay13}} (Group-Weights). Let $W_{1}$ and $W_{2}$ be admissible polynomials which have the same weights. Suppose $G \le G_{W_{1}}^{max}$ and $G \le G_{W_{2}}^{max}$. Then $\A_{W_{1},G} \isom \A_{W_{2},G}$.

\vspace{0.1 in}

\noindent Note that one can give the $\A$-model a product and pairing such that $\A$ is a Frobenius algebra. The above is then an isomorphism of Frobenius algebras, not just graded vector spaces.

\subsection{Properties of Invertible Polynomials}

\noindent Our initial intuition tells us that some of the properties of invertible polynomials should extend to the noninvertible case. For example, we'd like to keep the results of the following proposition.
\vspace{0.1 in}

\noindent \textbf{Proposition}. Let $W$ be an invertible polynomial. Then

(1) $W$ and $W^{T}$ have the same number of variables.

(2) $\left( G_{W}^{max} \right)^{T} = \{0\}$. 

(3) $\A_{W, G_{W}^{max}} \isom \B_{W^{T},\{0\}}$, as graded vector spaces.

\vspace{0.1 in}

\noindent \textbf{Proof}. (1) follows from noticing that the exponent matrix of $W$ is  square. Hence its transpose is also square and of the same size, so $W$ and $W^{T}$ have the same number of variables. (2) was stated previously in Section 2.3. (3) is a special case of the mirror symmetry conjecture that has been verified. Reference Theorem 4.1 in \cite{Kra09}. $\Box$
\vspace{0.1 in}

\noindent Part (3) of the proposition is especially important, and will be what we use to look for candidate transpose polynomials. In other words, given a noninvertible polynomial $W$, we would like to identify a candidate polynomial $W^{T}$ that satisfies $\mathlarger{\bigoplus\limits_{g \in G_{W}^{max}}} \left( \mathcal{Q}_{W|_{\text{fix}(g)}} \right)^{G_{W}^{max}} \isom \mathcal{Q}_{W^{T}}$. Though we would like this isomorphism to hold for all levels of algebraic structure, we will mainly investigate it on the level of graded vector spaces. For the benefit of the reader, we will restate the first conjecture.
\vspace{0.1 in}

\noindent \textbf{Conjecture 1}. For any admissible polynomial $W$ in $n$ variables, there exists a corresponding admissible polynomial $W^{T}$ in $n$ variables satisfying $\A_{W,G_{W}^{max}} \isom \B_{W^{T}, \{0\}}$.

\section{Results}

\subsection{Disproving Conjecture 1}

To disprove Conjecture 1, we prove a related nonexistence result. Note that this theorem is about any $W,\langle J \rangle$, whereas Conjecture 1 is about $W, G_{W}^{max}$. 

\vspace{0.1 in}

\noindent \textbf{Theorem}.  For any $n \in \N$, $n > 3$, let $W$ be an admissible but noninvertible polynomial in two variables with weight system $J = \left( \frac{1}{n},\;\frac{1}{n}\right)$, and let $G = \langle J \rangle$. Then there does not exist a corresponding $W^{T}$ in two variables satisfying $\mathcal{A}_{W,G} \isom \mathcal{B}_{W^{T}, \{0\} }$.
\vspace{0.1 in}

\noindent Before proving this theorem, we will demonstrate the hypothesis by exhibiting a few examples of such admissible polynomials for small values of $n$.
\begingroup
\renewcommand{\arraystretch}{1.2}%
\begin{align*}
\begin{array}{|l|c|l|l|c|l|}
\hline
n & J & \text{Some Examples} & n & J & \text{Some Examples} \\ \hline
 &  & x^{4} + y^{4} + x^{3}y &  &  & x^{5} + y^{5} + x^{4}y \\
4 & \left( \frac{1}{4},\;\frac{1}{4}\right) & x^{4} + x^{2}y^{2} + xy^{3} & 5 & \left( \frac{1}{5},\;\frac{1}{5}\right) &  x^{4}y + xy^{4} + x^{3}y^{2} + x^{2}y^{3}\\
 & &  x^{4} + xy^{3}  & &  & x^{5} + x^{2}y^{3} + xy^{4}\\ \hline
 &  & x^{6} + y^{6} + x^{5}y & & & x^{7} + y^{7} + x^{6}y \\
6 & \left( \frac{1}{6},\;\frac{1}{6}\right) & x^{5}y + x^{4}y^{2} + y^{6} & 7 & \left( \frac{1}{7},\;\frac{1}{7}\right) & x^{6}y + x^{5}y^{2} + y^{7}\\
  & &x^{6} + x^{2}y^{4} + xy^{5} + y^{6}  & & & x^{6}y + xy^{6}\\
\hline
\end{array}
\end{align*}
\endgroup

\vspace{0.1 in}

\noindent \textbf{Proof}. The idea of this proof is to choose an admissible polynomial with weight system $J = \left( \frac{1}{n},\;\frac{1}{n}\right)$, compute some formulas for its $\A$-model using the group $\langle J \rangle$, and show that there is no corresponding isomorphic unorbifolded $\B$-model. Then, under the Group-Weights isomorphism for $\A$-models, we will be able to generalize the result for any admissible polynomial with the same weights.  

\vspace{0.1 in}

\noindent To start,  we need an admissible polynomial in two variables with weight system $J = \left( \frac{1}{n},\;\frac{1}{n}\right)$. Let $W' = x^{n} + y^{n} + x^{n-1}y$, and let $G = \langle J \rangle$. Certainly $W'$ has weight system $J$, and $G$ fixes $W'$. 
\vspace{0.1 in}

\noindent For the unorbifolded $\B$-model, we know that $\dim(\B_{W^{T},\{0\} }) = \prod_{i = 1}^{n} \left( \frac{1}{q_{i}} - 1 \right)$ and that the highest degree of its graded pieces is given by $2\sum_{i = 1}^{n} \left(1 - 2q_{i} \right)$. In order to have $\A_{W,G} \isom \B_{W^{T}, \{0\} }$, we need the degrees of the vector spaces and the degrees of each of the graded pieces to be equal. Therefore we now need corresponding formulas for the dimension of the $\A$-model vector space and the degree of the highest degree piece of the $\A$-model. 
\vspace{0.1 in}

\noindent \textbf{Lemma}. As a graded vector space, $\dim\left(\A_{W',G} \right) = 2n - 2$, and the highest degree of any element is $\frac{2(2n-4)}{n}$. ($n \in \N$, $n \ge 3$). 
\vspace{0.1 in}

\noindent \textbf{Proof of Lemma}. Recall that $\A_{W',G} = \mathlarger{\bigoplus\limits_{g \in G}} \left( \mathcal{Q}_{W'|_{\text{fix}(g)}} \right)^{G}$. Notice that in our case $G = \langle \left( \frac{1}{n},\frac{1}{n} \right)  \rangle = \{(0,0)$, $\left( \frac{1}{n}, \frac{1}{n} \right)$, $\dots$, $\left( \frac{n-1}{n}, \frac{n-1}{n} \right)  \}$. Then $W'|_{\text{fix}(g)} = W'$ only for $g = (0,0)$. Otherwise $W'|_{\text{fix}(g)}$ is trivial.

\begin{description}
\item[Case 1] When $W'|_{\text{fix}(g)}$ is trivial, we get $n-1$ basis elements of the form $\lfloor 1; g \rceil$.

\item[Case 2] $W'|_{\text{fix}(g)} = W'$. Then $g = (0,0)$. The basis elements we get in this case are of the form $\lfloor x^{a}y^{b}; (0,0)\rceil$ where $a + b \equiv n - 2 \mod n$ and $a, b \in \{0, 1, \dots, n-2\}$. So we have $(a,b) = (0, n-2), (1,n-3), \dots, (n-3, 1), (n-2,0)$. Hence there are $n-1$ basis elements of this type.
\end{description}
The total dimension of $\A_{W',G}$ is therefore $(n - 1) + (n - 1) = 2n-2$. 

\vspace{0.1 in}

\noindent Now we will consider the degree of each basis element. Recall that 
\begin{align*}
\deg( \lfloor m ; g \rceil ) =  \dim(\text{fix}(g)) + 2 \sum_{i=1}^{n}(g_{i} - q_{i}),
\end{align*}
where $g = (g_{1},\dots,g_{n})$ and $J = (q_{1},\dots,q_{n})$ is the vector of quasihomogeneous weights.
\vspace{0.1 in}

\noindent For $g = (0,0)$, the degree is $ 2 + \left(-\frac{2}{n}\right) + \left(-\frac{2}{n}\right) = \frac{2(n-2)}{n}$. Also notice by the above equation that $\deg \left( \lfloor 1; \left( \frac{n-1}{n}, \frac{n-1}{n} \right) \rceil \right) > \deg \left( \lfloor 1; \left( \frac{m}{n}, \frac{m}{n} \right) \rceil\right)$ for all $m \in \{1, \dots, n-2\}$. Compute $\deg \left( \lfloor 1; \left( \frac{n-1}{n}, \frac{n-1}{n} \right) \rceil \right)= \frac{2(2n-4)}{n}$, and notice that $\frac{2(2n-4)}{n} = 2\left(\frac{2(n-2)}{n} \right) > \frac{2(n-2)}{n}$ for all $n \ge 3$. Hence the degree of the highest degree part of $\A_{W,G}$ is $\frac{2(2n-4)}{n}$. $\Box$
\vspace{0.1 in}

\noindent From the lemma, we now have the following system of equations for the possible weights $q_{1}, q_{2}$ for a candidate $W^{T}$:
\begin{align*}
\left(\frac{1}{q_{1}} - 1 \right) \left(\frac{1}{q_{2}} - 1 \right) &= 2n-2,\\
2\left((1 - 2q_{1}) + (1-2q_{2})\right) &= \frac{2(2n-4)}{n}.
\end{align*}
Solving for $q_{1}$ in the second equation, we have $q_{1} = \frac{2}{n} - q_{2}$. Substituting back into the first equation yields
\begin{align*}
n(2n-3)q_{2}^{2} + 2(3-2n)q_{2} + n - 2 = 0. 
\end{align*}
We now have a quadratic equation in $q_{2}$. Consider the discriminant 
\begin{align*}
D = -4(2n^{3} - 11n^{2} + 18n - 9).
\end{align*}
When $D < 0$, we will not have a real-valued solution for $q_{2}$. The above equation is a cubic polynomial that has roots at $n = 1, \frac{3}{2}, 3$. Since $D < 0$ for all $n > 3$, $q_{2}$ will not be real-valued for all $n > 3$. Thus there are no rational-valued solutions for the quasihomogeneous weights in this case. 
\vspace{0.1 in}

\noindent This shows that there is no $W^{T}$ in two variables satisfying $\A_{W',G} \isom \B_{W^{T}, \{0\}}$. Extending by the Group-Weights theorem, for any admissible polynomial $W$ with weights $\left( \frac{1}{n},\frac{1}{n} \right)$, we have that $\A_{W,G} \isom \A_{W',G}$. By this isomorphism, we know that $\dim\left(\A_{W,G } \right) = 2n - 2$ and the degree of its highest sector is $\frac{2(2n-4)}{n}$. Therefore, by what we have just shown, there cannot not exist any $W^{T}$ in two variables such that $\A_{W,G} \isom \B_{W^{T}, \{0\}}$. This proves the theorem. $\Box$
\vspace{0.1 in}

\noindent We do have the following solutions for $n \in \{1,2,3\}$. $n = 1$ yields the solution $\mathbf{q} = (1,1)$, $n= 2$ yields solutions $\mathbf{q} = (1,0)$, $(0,1)$, and $n=3$ gives a solution $\mathbf{q} = \left( \frac{1}{3}, \frac{1}{3} \right)$. However, since each coordinate must be in the interval (0, 1/2], $\mathbf{q} = \left( \frac{1}{3}, \frac{1}{3} \right)$ is the only valid weight system.
\vspace{0.1 in}

\noindent Our original conjecture (Conjecture 1) about the transpose of a noninvertible polynomial was that $W$ and $W^{T}$ have the same number of variables and $\left(G_{W}^{max} \right)^{T} = \{0\}$. We will now state a corollary to demonstrate that one of these assumptions must be false. 

\vspace{0.1 in}

\noindent \textbf{Corollary}. For any $n \in \N$, $n > 3$, let $W$ be a noninvertible polynomial in two variables with weight system $J = \left( \frac{1}{n},\;\frac{1}{n}\right)$ and $G_{W}^{max} = \langle J \rangle$. Then there does not exist a corresponding $W^{T}$ in two variables satisfying $\mathcal{A}_{W,G_{W}^{max}} \isom \mathcal{B}_{W^{T}, \{0\} }$.
\vspace{0.1 in}

\noindent The proof follows from the fact that for $W' = x^{n} + y^{n} + x^{n-1}y$ we have $\langle J \rangle = G_{W'}^{max}$. 
\vspace{0.1 in}

\noindent \textbf{Lemma}. The polynomial $W'$ has $G_{W'}^{max} = \langle J \rangle = \langle \left( \frac{1}{n},\;\frac{1}{n}\right) \rangle$ for all $n \in \N$, $n \ge 3$. 
\vspace{0.1 in}

\noindent The proof of the lemma relies on a theorem due to Lisa Bendall. We will state Bendall's theorem here, and refer the reader to the Appendix for a proof. 

\vspace{0.1 in}

\noindent \textbf{Theorem}. Let $W = x^{p} + y^{q}$. If a monomial satisfying the quasihomogeneous weights of $W$ is added to make the new polynomial $W' = x^{p} + y^{q} + x^{r}y^{s}$, then $G_{W'}^{max} = \langle (1/p,1/q),(1/n,0)\rangle$, where $n = \gcd(p,r)$. Alternatively, $G_{W'}^{max} = \langle (1/p,1/q),(0,1/m)\rangle $, where $m = \gcd(q,s)$. 

\vspace{0.1 in}

\noindent \textbf{Proof of Lemma}. By Lisa Bendall's theorem (see Appendix), $G_{W'}^{max} = \langle \left( \frac{1}{n},\;\frac{1}{n}\right), \left(0, \frac{1}{\gcd(n,1)} \right) \rangle =  \langle \left( \frac{1}{n},\;\frac{1}{n}\right) \rangle$, since $\gcd(n,1)$ = 1 and the generator $(0,1) \equiv (0,0) \mod 1$ contributes nothing. $\Box$
\vspace{0.1 in}

\noindent Since $W'$ has $G_{W'}^{max} = \langle J \rangle$, and since $W'$ satisfies the hypotheses of the previous theorem, we conclude that there does not exist a corresponding $W^{T}$ in two variables satisfying $\A_{W',G_{W'}^{max}} \isom \B_{W^{T},\{0\}}$. Extending by the Group-Weights theorem shows that any noninvertible $W$ with weights $J$ and $G_{W}^{max} = \langle J \rangle$ fails to have a $W^{T}$ in two variables satisfying the mirror symmetry alignment stated in the Corollary.



\subsection{Evidence Against Conjecture 2}

\noindent We will now consider finding a suitable $W^{T}$ in a different number of variables. By relaxing the constraint on the number of variables required in Conjecture 1, it is natural to make the following conjecture.

\vspace{0.1 in}

\noindent \textbf{Conjecture 2}. For any admissible $W$, there is a corresponding admissible $W^{T}$ satisfying $\A_{W,G_{W}^{max}} \isom \B_{W^{T}, \{0\}}$.

\vspace{0.1 in}

\noindent The following theorem is a start to disproving this conjecture.

\vspace{0.1 in}

\noindent \textbf{Theorem}.  For any admissible polynomial $W$ with weight system $J = \left( \frac{1}{5}, \frac{1}{5} \right)$ and $G = \langle J \rangle$, there is no corresponding admissible $W^{T}$ in 1, 2, or 3 variables satisfying $\A_{W,G} \isom \B_{W^{T}, \{0\}}$. 
\vspace{0.1 in}

\noindent \textbf{Proof}. For $W$ as given in the hypothesis, we have previously shown that the degree of the $\A$-model is 8, and the degree of its highest sector is $12/5$.
\vspace{0.1 in}

\noindent We will rule out the existence of a $W^{T}$ in these three cases. In one variable, we can only have $W^{T} = x^{9}$ to give us an unorbifolded $\B$-model of dimension 8. Then $q_{1} = \frac{1}{9}$, but $1 - \frac{2}{9} = \frac{7}{9} \ne \frac{6}{5}$. The two variable case is done by the previous theorem. 
\vspace{0.1 in}

\noindent  Now let $n \in \N$, $n \ge 3$. We have the following equations for a candidate weight system:
\begin{align}
\left( \frac{1}{q_{1}} - 1 \right) \left( \frac{1}{q_{2}} - 1 \right) \prod\limits_{i = 3}^{n} \left( \frac{1}{q_{i}} - 1 \right) &= 8, \\
2\left[(1 - 2q_{1}) + (1 - 2q_{2}) + \sum\limits_{i = 3}^{n} (1 - 2q_{i})\right] &= \frac{12}{5} .
\end{align}
Letting $A = 1 - \dfrac{8}{ \prod\limits_{i = 3}^{n} \left( \frac{1}{q_{i}} - 1 \right)   }$, and $B = \dfrac{5n - 6}{10} - \sum\limits_{i = 3}^{n} q_{i}$, equations (1) and (2) simplify to
\begin{align}
Aq_{1}q_{2} - q_{1} - q_{2} + 1 &= 0,\\
-q_{1} + B &= q_{2}.
\end{align}

\vspace{0.1 in}

\noindent For any $q_{i} \in (0, 1/2]$, we have that $\frac{1}{q_{i}} - 1 \ge 1$. By equation (1), we require that $\prod\limits_{i = 3}^{n} \left( \frac{1}{q_{i}} - 1 \right) \le 8$. This tells us that $1 \le \prod\limits_{i = 3}^{n} \left( \frac{1}{q_{i}} - 1 \right) \le 8$. Therefore we have that $-7 \le A \le 0$. 
\vspace{0.1 in}

\noindent From equation (2) we also have that $\sum\limits_{i = 3}^{n} (1 - 2q_{i}) \le \frac{6}{5}$. Rewriting the left-hand side gives us $(n-2) - 2\sum\limits_{i = 3}^{n} q_{i} \le \frac{6}{5}$. Subtracting $n-2$ from both sides yields $-\sum\limits_{i = 3}^{n} q_{i} \le \frac{16-5n}{10}$.


\vspace{0.1 in}

\noindent Substituting this into $B$ gives us
\begin{align*}
B = \dfrac{5n - 6}{10} - \sum\limits_{i = 3}^{n} q_{i} \le \dfrac{5n - 6}{10} + \dfrac{16-5n}{10} = 1. 
\end{align*}

\vspace{0.1 in}

\noindent Though we have developed the previous formulas in general, we will now restrict our attention to the case $n = 3$. When $A \ne 0$, we can use the quadratic formula to plot the real-valued solutions of $q_{1}$. In three variables, the discriminant $D = (AB)^{2} - 4A(B-1) \ge 0$ for $q_{3} \le 1/9$. This yields the following:

\begin{figure}[H]
\centering
\begin{minipage}{.5\textwidth}
  \centering
  \includegraphics[width=\linewidth]{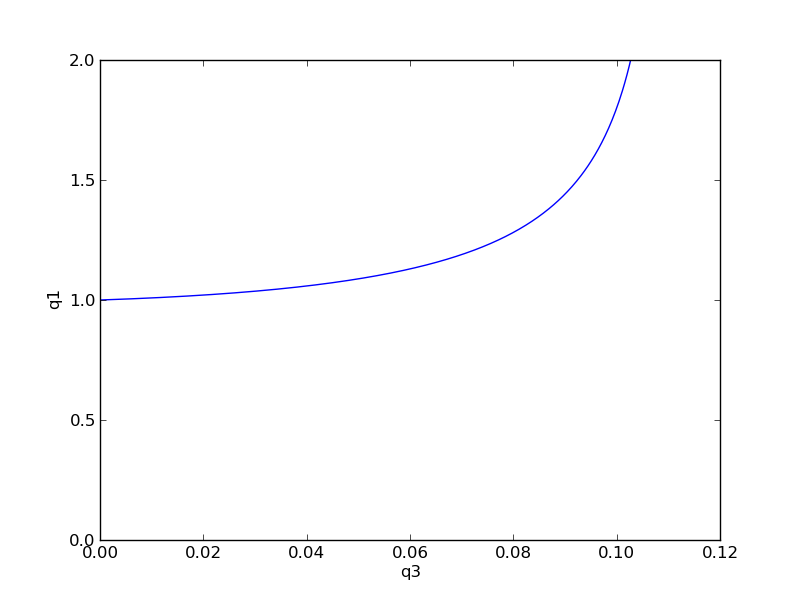}
\caption{\label{fig:pos}Positive solutions for $q_{1}$ \hspace{\textwidth}in the quadratic system (3) and (4)}
\end{minipage}%
\begin{minipage}{.5\textwidth}
  \centering
  \includegraphics[width=\linewidth]{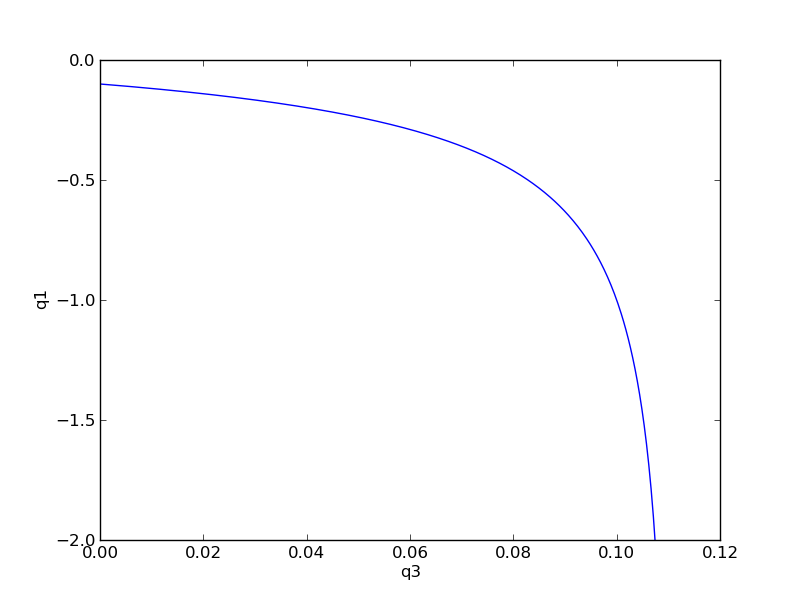}
\caption{\label{fig:neg}Negative solutions for $q_{1}$ \hspace{\textwidth}in the quadratic system (3) and (4)}
\end{minipage}
\end{figure}

\noindent None of these values of $q_{1}$ is in the interval $(0, 1/2]$, let alone $(0, 1/2] \cap \Q$.

\vspace{0.1 in}
\noindent Now when $A = 0$, we must have that $ \frac{1}{q_{i}} - 1 = 8$. Therefore by equation (1) we can only have $q_{1} = q_{2} = 1/2$. But equations (1) and (2) show that if this is the case, then we could have found a satisfactory weight system in just 1 variable without considering $q_{1}$ and $q_{2}$. Since we have already ruled out the case $n = 1$, we conclude that there are no valid weight systems for $W^{T}$ in three variables. $\Box$

\vspace{0.1 in}

\noindent The previous result casts doubt on the validity of Conjecture 2. Using the formulas developed in the last theorem may be useful in proving the following statement. 

\vspace{0.1 in}

\noindent \textbf{Conjecture 3}.  For any admissible polynomial $W$ with weight system $J = \left( \frac{1}{5}, \frac{1}{5} \right)$ and $G = \langle J \rangle$, there is no corresponding admissible $W^{T}$ satisfying $\A_{W,G} \isom \B_{W^{T}, \{0\}}$. 

\vspace{0.1 in}

\noindent Proving Conjecture 3 will demonstrate that the mirror symmetry construction $\A_{W,G_{W}^{max}} \isom \B_{W^{T}, \{0\}}$ does not, in general, extend to noninvertible $W$.

\section{Conclusion}

Given a polynomial $W$ fixed by a weight system $J = \left(\frac{1}{n},\frac{1}{n}\right)$ and group $G = \langle J \rangle$, and $m \in \N$ representing the number of variables in a candidate $W^{T}$, it is impossible to construct $\A_{W,G} \isom \B_{W^{T},\{0\}}$ in the following cases:
\begin{align*}
\begin{array}{cc|cccc}
  &   &   & m &   &  \\ 
  &  & 1 & 2 & 3 & \dots \\
\hline
  & 4 &   & X  &   &       \\
n & 5 & X  & X  & X  &       \\   
  & 6 &   & X  &   &       \\ 
  & \vdots &   & X  &   &       
\end{array}
\end{align*}
\vspace{0.1 in}



\noindent These results show that our original intuition about invertible polynomials and their transposes does not extend well to the noninvertible case. Even at the level of graded vector spaces, simply allowing an invertible polynomial to have one extra monomial seems to break this mirror symmetry construction.

\vspace{0.1 in}

\noindent Though we have not completely ruled out the possibility of noninvertible polynomials having a transpose, we have shown that this problem is difficult and will require further research to fully elucidate it. 


\section{References}

\begingroup
\renewcommand{\section}[2]{}%

\endgroup


\section{Appendix}

\noindent We used the following result when proving the corollary in Section 3.1. The theorem and proof are due to Lisa Bendall. We will reproduce the entire proof here because it is not publicly available elsewhere.
\vspace{0.1 in}

\noindent \textbf{Theorem}. Let $W = x^{p} + y^{q}$. If a monomial satisfying the quasihomogeneous weights of $W$ is added to make the new polynomial $W' = x^{p} + y^{q} + x^{r}y^{s}$, then $G_{W'}^{max} = \langle (1/p,1/q),(1/n,0)\rangle$, where $n = \gcd(p,r)$. Alternatively, $G_{W'}^{max} = \langle (1/p,1/q),(0,1/m)\rangle $, where $m = \gcd(q,s)$. 
\vspace{0.1 in}

\noindent \textbf{Proof}. Any element $(\theta_{1}, \theta_{2}) \in G_{W'}^{max}$ must satisfy the following matrix equation:
\begin{align*}
\begin{bmatrix} p&0\\0&q \\r&s \end{bmatrix} \begin{bmatrix} \theta_{1} \\ \theta_{2} \end{bmatrix} = \begin{bmatrix}k_{1} \\ k_{2} \\ k_{3} \end{bmatrix} \in \Z^{3}.
\end{align*}
This yields the following three equations:
\begin{align*}
\theta_{1} = \frac{k_{1}}{p}, \hspace{50 pt} \theta_{2} = \frac{k_{2}}{1}, \hspace{50 pt} r\theta_{1} + s\theta_{2} = k_{3}.
\end{align*}

\noindent We also note that since the new monomial satisfies the weights vector of $W$, then $\frac{r}{p} + \frac{s}{q} = 1$. From this equation, it follows that $s = \frac{pq - rq}{p}$, which we can substitute along with the first two equations into the last equation to get the equation $r(k_{1} - k_{2}) = p(k_{3} - k_{2})$. Dividing out by the gcd of $r$ and $p$, we get the equation $r'(k_{1} - k_{2}) = p'(k_{3} - k_{2})$, where $r'$ and $p'$ are relatively prime.
\vspace{0.1 in}

\noindent From this, we know that $p' \mid (k_{1} - k_{2})$, or in other words, $k_{1} - k_{2} = k_{4} p'$ for some $k_{4} \in \Z$. Now, dividing both sides by $p$, we get $\frac{k_{1}}{p} - \frac{k_{2}}{p} = \frac{k_{4}}{n}$. Next, substitute $p\theta_{1}$ for $k_{1}$. We find that $\theta_{1} = k_{2}\left(\frac{1}{p}\right) + k_{4} \left(\frac{1}{n}\right)$. From the second equation, we already know that $\theta_{2} = \frac{k_{2}}{q}$, so we have the following equation in vector form:
\begin{align*}
(\theta_{1},\theta_{2}) = k_{2}(1/p, 1/q) + k_{4}(1/n, 0),
\end{align*}
where $k_{2}, k_{4} \in \Z$.
\vspace{0.1 in}

\noindent Now, to show that this generates the group, we show that anything of form $k_{2}(1/p,1/q) + k_{4}(1/n,0)$ satisfies the three original equations for some three arbitrary integers. For the first equation, note that $1/n = p'/p$, thus $\theta_{1} = (k_{2} + k_{4} p')(1/p)$, so it is satisfied for some integer. The second equation follows immediately. For the final equation, plugging in we get $r(k_{2}/p + k_{4}/n) + s(k_{2}/q) = (r/p + s/q)k_{2} + (r/n)k_{4} = k_{2} + r'k_{4} \in \Z$. Therefore, any element of the form $k_{2}(1/p,1/q) + k_{4}(1/n,0)$ is in $G_{W'}^{max}$. Thus $G_{W'}^{max} = \langle (1/p,1/q),(1/n,0)\rangle$.
\vspace{0.1 in}

\noindent Note:  by substituting $r = \frac{pq - sp}{q}$ in the last equation, we  get the alternate set of generators $(1/p,1/q)$ and $(0,1/m)$. $\Box$

%
%
%
%
%

\end{document}